\documentclass[11pt,a4paper,leqno]{article}
\usepackage[latin2]{inputenc}
\usepackage{latexsym}
\usepackage{amssymb,amsmath,amsopn}
\usepackage{mathrsfs}
\usepackage[dvips]{graphicx}   
\usepackage{color,epsfig}      

\newtheorem{lem}{\it LEMMA}

\newtheorem{thm}{\it THEOREM}
\newtheorem{rem}{\it REMARK}

\newcommand{\B}{{\mathscr B}}
\newcommand{\C}{{\mathscr C}}

\begin{document}


\centerline{\Large  \bf PERFECT PACKING OF SQUARES}

\vspace{10pt}

\centerline{ANTAL JO\' OS}


\begin{abstract}
{\small It is known that $\sum\limits_{i =1}^\infty {1/
i^2}={\pi^2/6}$. Meir and Moser asked what is the smallest
$\epsilon$ such that all the squares of sides of length $1$, $1/2$,
$1/3$, $\ldots$ can be packed into a rectangle of area
${\pi^2/6}+\epsilon$. A packing into a rectangle of the right area
is called perfect packing. Chalcraft packed the squares of sides of
length $1$, $2^{-t}$, $3^{-t}$, $\ldots$ and he found perfect
packing for $1/2<t\le3/5$. We will show based on an algorithm by
Chalcraft that there are perfect packings if $1/2<t\le2/3$. Moreover
we show that there is a perfect packing for all $t$ in the range
$\log_32\le t\le2/3$.}
\end{abstract}




\medskip

\centerline {\bf 1. INTRODUCTION}

\medskip

Meir and Moser \cite{MeirMoser1968} originally noted that since
$\sum\limits_{i =2}^\infty {1/ i^2}={\pi^2/6}-1$, it is reasonable
to ask whether the set of squares with sides of length $1/2$, $1/3$,
$1/4$, $\ldots$ can be packed into a rectangle of area
${\pi^2/6}-1$. Failing that, find the smallest $\epsilon$ such that
the squares can be packed in a rectangle of area
${\pi^2/6}-1+\epsilon$. The problem also appears in \cite{Croft},
\cite{BMP}, \cite{BojuFunar}.

A packing into a rectangle of the right (resp. not the right) area
is called \emph{perfect} (resp. \emph{imperfect}) packing. In
\cite{MeirMoser1968}, \cite{Jennings}, \cite{Ball}, \cite{Paulhus}
can be found better and better imperfect packing.

Chalcraft \cite{Chalcraft} generalized this question. He packed the
squares of side $n^{-t}$ for $n=1,2,\ldots$ into a square of the
right area. He proved that for all $t$ in the range $[0.5964,0.6]$
there is a perfect packing of the squares. In \cite{Chalcraft} can
be read that "Other packings will work for other ranges of $t$. We
can probably make the $t_0$ in Theorem 8 as close to $1/2$ as
desired in this way. The more interesting challenge, however, seems
to be to increase the bound $t\le3/5$." Our aim is to increase this
bound.\\
W\"astlund \cite{Wastlund} proved if $1/2<t<2/3$, then the squares
of side $n^{-t}$ for $n=1,2,\ldots$ can be packed into some finite
collection of square boxes of the same area $\zeta(2t)$ as the total
area of the tiles. This is an increase of the bound $t\le3/5$, but
we have many enclosing rectangles.\\
We can find several papers in this topic e.g. \cite{Martin},
\cite{Altunkaynak}, \cite{Liu}.

\medskip

\medskip

\centerline {\bf 2. PERFECT PACKING}

\medskip

\begin{thm}\label{t1}
For $t=2/3$, the squares $S^t_n$ $(n\ge1)$ can be packed perfectly
into the rectangle of dimensions $\zeta(2t)\times 1$.
\end{thm}

\begin{thm}\label{t2}
For all $t$ in the range $\log_32\le t\le2/3$, the squares $S^t_n$
$(n\ge1)$ can be packed perfectly into the rectangle of dimensions
$\zeta(2t)\times 1$.
\end{thm}

\medskip

\medskip

\centerline {\bf 3. NOTATION}

\medskip

We use the Chalcraft$'$s algorithm in \cite{Chalcraft} and we modify
the proof of Chalcraft. For the sake of simplicity we use the
Chalcraft$'$s notation. For the completeness, we recall these.

Throughout the paper the width of a rectangle will always refer to
the shorter side and the height will always refer to the longer
side. We use the constant $1/2<t\le 2/3$. As usual,
$\zeta(t)=\sum\limits_{i=1}^\infty i^{-t}$. Let $S^t_n$ denote the
square of side length $n^{-t}$. A box is a rectangle of sides $x,
y>0$. Let $x\times y$ denote the box $B$ of sides $x$ and $y$. We
define its area $a(B)=xy$, its semi-perimeter $p(B)=x+y$, its width
$w(B)=\min(x, y)$ and its height $h(B)=\max(x, y)$. Given a set of
boxes $\B=\{B_1,\ldots, B_n \}$, we define
$a(\B)=\sum\limits_{i=1}^n a(B_i)$, $h(\B)= \sum\limits_{i=1}^n
h(B_i)$ and $w(\B)=\max\limits_{i=1,\ldots,n} w(B_i)$. Let
$a(\emptyset)= h(\emptyset)=w(\emptyset)=0$.

\medskip

\medskip

\centerline {\bf 4. CHALCRAFT'S ALGORITHM}

\medskip

For the completeness, we repeat the description of the Chalcraft's
algorithm.

First we recall the subroutine of Chalcraft, which we call Algorithm
$\bf b$ as in \cite{Chalcraft}.

\noindent \emph{Algorithm} $\bf b$

Input: An integer $n\ge1$ and a box $B$, where $w(B)=n^{-t}$.

Output: If the algorithm terminates, then it defines an integer
$m_{\bf b}\linebreak[4]=m_{\bf b} (n, B)>n$ and a set of boxes
$\B_{\bf b}=\B_{\bf b} (n, B)$.

Action: If the algorithm terminates, then it packs the squares
$S^t_n , \ldots, S^t_{m_{\bf b}-1}$ into $B$, and $\B_{\bf b}$ is
the set of boxes containing the remaining area. If it does not
terminate, then it packs the squares $S^t_n, S^t_{n+1}, \ldots$ into
$B$.

\noindent
\begin{tabular}{p{0.8cm}p{10.7cm}}
({\bf b}1)& Let $n_1=n+1$, $x_1=h(B)-n^{-t}$ and $\B_1=\emptyset$. \\
\end{tabular}

\noindent
\begin{tabular}{p{0.8cm}p{10.7cm}}
({\bf b}2)& Put the square $S^t_n$ snugly at one end of $B$. \\
\end{tabular}

\noindent
\begin{tabular}{p{0.8cm}p{10.7cm}}
({\bf b}3)& If $x_1>0$, then let $B_1$ be the remainder of $B$, so
that $B_1$ has dimensions $x_1\times n^{-t}$.\\
\end{tabular}

\noindent
\begin{tabular}{p{0.8cm}p{10.7cm}}
({\bf b}4)& For $i=1, 2, \ldots$ \\
\end{tabular}

\noindent
\begin{tabular}{p{1.2cm}p{10.2cm}}
({\bf b}5) &  (Note: At stage $i$, we have packed $S^t_n,\ldots,
S^t_{n_i-1}$ into $B$. The remaining boxes are $\B_i$, which we
never use again in this algorithm, and $B_i$ (as
long as $x_i>0$), which has dimensions $x_i\times n^{-t}$.)\\
\end{tabular}

\noindent
\begin{tabular}{p{1.2cm}p{10.2cm}}
({\bf b}6)&  If $x_i=0$, then terminate with $m_{\bf b}=n_i$
and $\B_{\bf b}=\B_i$.\\
\end{tabular}

\noindent
\begin{tabular}{p{1.2cm}p{10.2cm}}
({\bf b}7)&  If $x_i<n^{-t}_i$, then terminate with $m_{\bf
b}=n_i$ and $\B_{\bf b}=\B_i\cup\{B_i\}$.\\
\end{tabular}

\noindent
\begin{tabular}{p{1.2cm}p{10.2cm}}
({\bf b}8)&  Let $x_{i+1}=x_i-n^{-t}_i$.\\
\end{tabular}

\noindent
\begin{tabular}{p{1.2cm}p{10.2cm}}
({\bf b}9) &  If $x_{i+1}=0$, then let $C_i=B_i$.\\
\end{tabular}

\noindent
\begin{tabular}{p{1.2cm}p{10.2cm}}
({\bf b}10)& If $x_{i+1}>0$, then split $B_i$ into two boxes: one
called $C_i$ with dimensions $n^{-t}_i\times n^{-t}$, and the
other called $B_{i+1}$ with dimensions  ${x_{i+1}}\times n^{-t}$.\\
\end{tabular}

\noindent
\begin{tabular}{p{1.2cm}p{10.2cm}}
({\bf b}11)&  Apply Algorithm $\bf b$ recursively with inputs $n_i$
and $C_i$. If this terminates, let $n_{i+1}=m_{\bf b}(n_i ,
C_i)$ and $\C_i=\B_{\bf b} (n_i , C_i)$.\\
\end{tabular}

\noindent
\begin{tabular}{p{1.2cm}p{10.2cm}}
({\bf b}12)&  Let $\B_{i+1}=\B_i\cup \C_i$.\\
\end{tabular}

\noindent
\begin{tabular}{p{0.8cm}p{10.7cm}}
({\bf b}13)&  End For.\\
\end{tabular}

The subroutine $\bf b$ is used in the Chalcraft's algorithm $\bf c$.

\noindent \emph{Algorithm} $\bf c$

Input: An integer $n\ge1$ and a set of boxes $\B$.

Action: If the algorithm does not fail, then it packs the squares
$S^t_n, S^t_{n+1}, \ldots$ into $\B$.

\noindent
\begin{tabular}{p{0.8cm}p{10.7cm}}
({\bf c}1)& Let $n_1=n+1$ and $\B_1=\B$. \\
\end{tabular}

\noindent
\begin{tabular}{p{0.8cm}p{10.7cm}}
({\bf c}2)& For $i=1, 2, \ldots$ \\
\end{tabular}

\noindent
\begin{tabular}{p{1.2cm}p{10.2cm}}
({\bf c}3)& (Note: At stage $i$, we have packed $S^t_n,\ldots,
S^t_{n_i-1}$ into $B$. The remaining boxes are $\B_i$.)\\
\end{tabular}

\noindent
\begin{tabular}{p{1.2cm}p{10.2cm}}
({\bf c}4)& If $w(\B_i)<n^{-t}_i$, then fail.\\
\end{tabular}

\noindent
\begin{tabular}{p{1.2cm}p{10.2cm}}
({\bf c}5)& Let $w_i=\min\{w(C) | C\in \B_i , w(C)\ge n^{-t}_i
\}$.\\
\end{tabular}

\noindent
\begin{tabular}{p{1.2cm}p{10.2cm}}
({\bf c}6)& Let $h_i=\min\{h(C) | C\in \B_i , w(C)=w_i
\}$.\\
\end{tabular}

\noindent
\begin{tabular}{p{1.2cm}p{10.2cm}}
({\bf c}7)& Choose any $B_i\in \B_i$ which satisfies $w(B_i)=w_i$
and $h(B_i)=h_i$.\\
\end{tabular}

\noindent
\begin{tabular}{p{1.2cm}p{10.2cm}}
({\bf c}8)& If $w_i=h_i=n^{-t}_i$ , then\\
\end{tabular}

\noindent
\begin{tabular}{p{1.6cm}p{9.8cm}}
({\bf c}9)& Put $S^t_{n_i}$ snugly into $B_i$.\\
\end{tabular}

\noindent
\begin{tabular}{p{1.6cm}p{9.8cm}}
({\bf c}10)& Let $\B_{i+1}=\B_i\setminus \{B_i\}$.\\
\end{tabular}

\noindent
\begin{tabular}{p{1.6cm}p{9.8cm}}
({\bf c}11)& Let $n_{i+1}=n_i+1$.\\
\end{tabular}

\noindent
\begin{tabular}{p{1.2cm}p{10.2cm}}
({\bf c}12)&  Else\\
\end{tabular}

\noindent
\begin{tabular}{p{1.6cm}p{9.8cm}}
({\bf c}13)& Cut $B_i$ into two boxes: one called $C_i$ of
dimensions $w_i\times n^{-t}_i$ and tho other
called $D_i$ of dimensions $w_i\times (h_i-n^{-t}_i)$.\\
\end{tabular}

\noindent
\begin{tabular}{p{1.6cm}p{9.8cm}}
({\bf c}14)& Call Algorithm $\bf b$ with inputs $n_i$ and $C_i$. If
this terminates, then
let $n_{i+1}=m_{\bf b}(n_i , C_i)$ and $\C_i=\B_{\bf b} (n_i , C_i)$.\\
\end{tabular}

\noindent
\begin{tabular}{p{1.6cm}p{9.8cm}}
({\bf c}15)& Let $\B_{i+1}=\B_i\setminus \{B_i\}\cup
\C_i\cup \{D_i\}$.\\
\end{tabular}

\noindent
\begin{tabular}{p{1.2cm}p{10.2cm}}
({\bf c}16)&  End If.\\
\end{tabular}

\noindent
\begin{tabular}{p{0.8cm}p{10.7cm}}
({\bf c}17)&  End For.\\
\end{tabular}

\medskip

\medskip

\centerline {\bf 5. THE PROOF}

\medskip

The key lemma of Chalcraft is Lemma 1 in \cite{Chalcraft}. We modify
that in the following way.

\begin{lem}
If $\B=\{B_1,\ldots,B_n\}$ $(n\ge1)$, then $a(\B)\le w(\B)h(\B)$.
\end{lem}

\noindent {\it Proof.}  We have
$$a(\B)=\sum\limits_{i=1}^n a(B_i)=\sum\limits_{i=1}^n w(B_i)h(B_i)\le \sum\limits_{i=1}^n
w(\B)h(B_i) $$$$= w(\B)\sum\limits_{i=1}^n h(B_i) =w(\B) h(\B),$$
which completes the proof.

We prove the modified Chalcraft's lemmas in which we use the height
instead of the semi-perimeter.

\begin{lem} Suppose $w(B)=n^{-t}$ and Algorithm $\bf b$ with inputs $n$ and $B$ terminates
with $m_{\bf b}=m_{\bf b} (n, B)$ and $\B_{\bf b}=B_{\bf b} (n, B)$.
Therefore $$h(\B_{\bf b})\le \sum\limits^{m_{\bf b}-1}_{j=n}
j^{-t}.$$
\end{lem}

\noindent {\it Proof.} The proof is similar to the proof of Lemma 2
in \cite{Chalcraft}. For completeness, we write it again.

The proof is by induction on the number of squares packed. Of
course, if $\bf b$ terminates with $m_{\bf b}=n+1$, then $h(\B_{\bf
b})\le n^{-t}= \sum\limits^{m_{\bf b}-1}_{j=n} j^{-t}.$\\
We can assume that the lemma is true of all the recursive calls to
Algorithm $\bf b$. We can also assume that $\bf b$ and all the
recursive calls to $\bf b$ terminated. Suppose Algorithm $\bf b$
terminates when $i=k$, so $m_{\bf b}=n_k$. Since Algorithm $\bf b$
terminated without placing the next square, $x_k<n^{-t}_k <n^{-t}$,
so $h(B_k)=n^{-t}$. Now by induction,
$$h(\C_i)\le \sum\limits^{n_{i+1}-1}_{j=n_i}
j^{-t} \quad  \hbox{ for } i<k,$$
$$\sum\limits^{k-1}_{i=1}h(\C_i)\le \sum\limits^{n_{k}-1}_{j=n_1}
j^{-t}= \sum\limits^{m_{\bf b}-1}_{j=n+1} j^{-t}.$$ If the condition
in (${\bf b}6$) was true, then
$$h(\B_{\bf b})=\sum\limits^{k-1}_{i=1} h(\C_i)
\le\sum\limits^{m_{\bf b}-1}_{j=n+1} j^{-t}<\sum\limits^{m_{\bf
b}-1}_{j=n} j^{-t}.$$ If the condition in (${\bf b}7$) was true,
then
$$h(\B_{\bf b})=\sum\limits^{k-1}_{i=1} h(\C_i)+h(B_k)
\le\sum\limits^{m_{\bf b}-1}_{j=n+1}
j^{-t}+n^{-t}=\sum\limits^{m_{\bf b}-1}_{j=n} j^{-t},$$ which
completes the proof.

\begin{lem} We have
\begin{equation} \label{eq:1}
(b+1)^{1-t}-a^{1-t}<(1-t)\sum\limits_{j=a}^b
j^{-t}<b^{1-t}-(a-1)^{1-t},
\end{equation}
\begin{equation} \label{eq:2}
a^{1-2t}-(b+1)^{1-2t}<(2t-1)\sum\limits_{j=a}^b
j^{-2t}<(a-1)^{1-2t}-b^{1-2t}.
\end{equation}
\end{lem}

\noindent {\it Proof.} We omit the proof.

\begin{lem} Consider step $({\bf c}4)$ for some value of $i$. Suppose the following
conditions hold.
\begin{equation} \label{eq:3}a(\B_i)\ge\sum\limits_{j=n_i}^\infty
j^{-2t}, \end{equation}
\begin{equation}\label{eq:4} h(\B_i)\le
{n^{1-t}_i\over 2t-1}. \end{equation} Therefore step $({\bf c}4)$
will not fail for this value of $i$.
\end{lem}

\noindent {\it Proof.} We assume, that the algorithm fail. Therefore
we have $w(\B_i)<n_i^{-t}$. By Lemma 1, (\ref{eq:4}), (\ref{eq:2}),
$$a(\B_i)\le w(\B_i)h(\B_i)<{n^{1-2t}_i\over 2t-1}\le
\sum\limits_{j=n_i}^\infty j^{-2t}\le a(\B_i),$$ a contradiction,
which completes the proof of the lemma.

\begin{lem} Given an integer $n\ge1$ and a non-empty set of boxes $\B$,
suppose the following conditions hold
\begin{equation} \label{eq:6}a(\B)\ge\sum\limits_{j=n}^\infty
j^{-2t}, \end{equation}
\begin{equation}\label{eq:7} h(\B)\le
{1\over 1-t}(n-1)^{1-t}, \end{equation}
$$t\le {2\over3}.$$
If we run Algorithm $\bf c$ with the inputs $n$ and $\B$, then the
conditions
\begin{equation} \label{eq:8}a(\B_i)\ge\sum\limits_{j=n_i}^\infty
j^{-2t}, \end{equation}
\begin{equation}\label{eq:9} h(\B_i)\le
h(\B)+\sum\limits_{j=n}^{n_i-1} j^{-t}. \end{equation} hold at step
$({\bf c}4)$ for all $i\ge1$ for which step $({\bf c}4)$ is
executed. Moreover, the algorithm will never fail.
\end{lem}

\noindent {\it Proof.} First, we will show that (\ref{eq:8}) and
(\ref{eq:9}) ensure that the algorithm will not fail. By
(\ref{eq:9}), (\ref{eq:1}), and (\ref{eq:7}),
$$h(\B_i)\le
h(\B)+\sum\limits_{j=n}^{n_i-1} j^{-t}$$
$$ < h(\B)+{1\over
1-t}((n_i-1)^{1-t}-(n-1)^{1-t})$$
$$\le{1\over1-t}(n_i-1)^{1-t}.$$
Since $t\le2/3$, $${1\over1-t}\le{1\over2t-1}.$$ Thus
$$h(\B_i)<{1\over1-t}(n_i-1)^{1-t}\le{1\over2t-1}(n_i-1)^{1-t}
<{n_i^{1-t}\over2t-1}.$$ By Lemma 4, $({\bf c}4)$ will not fail.

Now we prove (\ref{eq:8}) and (\ref{eq:9}) by induction on $i$. Of
course they hold for $i=1$ and (\ref{eq:8}) holds for all $i$. Let
$i>1$ be the smallest $i$ for which (\ref{eq:9}) is not true.

If the condition in $({\bf c}8)$ was true for $i-1$, then
$h(\B_i)=h(\B_{i-1})-n^{-t}_{i-1}$ and $n_i=n_{i-1}+1$. Thus by
induction,
$$h(\B_i)=h(\B_{i-1})-n^{-t}_{i-1}\le
h(\B)+\sum\limits_{j=n}^{n_{i-1}-1} j^{-t}-n^{-t}_{i-1}$$
$$<h(\B)+\sum\limits_{j=n}^{n_{i-1}-1}
j^{-t}=h(\B)+\sum\limits_{j=n}^{n_{i}-2}
j^{-t}<h(\B)+\sum\limits_{j=n}^{n_{i}-1} j^{-t}.$$ If the condition
in $({\bf c}8)$ was not true for $i-1$,
then we distinguish two cases.\\
If $w_{i-1}\ge
h_{i-1}-n^{-t}_{i-1}$ (that is $h(D_{i-1})=w_{i-1}$), then
$$h(\B_i)=h(\B_{i-1})+h(\C_{i-1})-h(B_{i-1})+h(D_{i-1})$$
$$=h(\B_{i-1})+h(\C_{i-1})-h_{i-1}+w_{i-1}\le h(\B_{i-1})+h(\C_{i-1}).$$
If $w_{i-1}< h_{i-1}-n^{-t}_{i-1}$ (that is
$h(D_{i-1})=h_{i-1}-n^{-t}_{i-1}$), then similarly
$$h(\B_i)
\le h(\B_{i-1})+h(\C_{i-1}).$$ By induction and Lemma 2,
$$h(\B_{i})\le h(\B_{i-1})+h(\C_{i-1})$$
$$\le h(\B)+\sum\limits_{j=n}^{n_{i-1}-1}
j^{-t}+\sum\limits_{j=n_{i-1}}^{n_{i}-1} j^{-t}
=h(\B)+\sum\limits_{j=n}^{n_{i}-1} j^{-t},
$$ which completes
the proof.

\begin{figure}
\center{
  \includegraphics[width=0.6\textwidth]{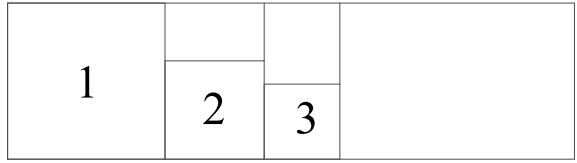}\\
  \caption{The squares $S_1^t,S_2^t,S_3^t$ and the set of boxes $\B$.}\label{fig10}}
\end{figure}

\noindent {\it Proof of Theorem \ref{t1}.} If the first three
squares are packed in the box $B=\zeta(2t)\times 1$ as in Fig.
\ref{fig10} (this is the Paulhus's algorithm \cite{Paulhus}),
 then the remaining boxes are
$$\B=\left\{ \left(\zeta(2t)-1-2^{-t}-3^{-t}\right)\times1, 2^{-t}\times \left(1-2^{-t}\right),
 3^{-t}\times \left(1-3^{-t}\right)\right\} $$ and
$$h(\B)=
\zeta(2t)-2\cdot3^{-t}=2.639$$$$< 4.327={1\over1-t}(4-1)^{1-t}.$$ By
Lemma 5, the Algorithm $\bf c$ pack perfectly the squares $S^t_n$
$(n\ge4)$ into $\B$, which completes the proof.

\begin{rem}
The squares $S^t_n$ $(n\ge1)$ in Theorem \ref{t1} can be packed
similarly in a square of the right area.
\end{rem}

\noindent {\it Proof of Theorem \ref{t2}.} If the first three
squares are packed in the box $B=\zeta(2t)\times 1$ as in Fig.
\ref{fig10}, then the remaining boxes are
$$\B=\left\{ \left(\zeta(2t)-1-2^{-t}-3^{-t}\right)\times1, 2^{-t}\times \left(1-2^{-t}\right),
 3^{-t}\times \left(1-3^{-t}\right)\right\}. $$ Observe $\zeta(2t)-1-2^{-t}-3^{-t}>1$,
$2^{-t}>1-2^{-t}$ and $1-3^{-t}\ge 3^{-t}$ if
$t\in\left[\log_32,2/3\right]$. Let $f(t)=h(\B)$. Thus
$$h(\B)=f(t)= \zeta(2t)-2\cdot3^{-t}.$$ Since
$$g(t)={1\over1-t}3^{1-t}$$ is an increasing, $f(t)$ is a
decreasing function on the interval $\left[\log_32,2/3\right]$ and
$$f(\log_32)=3.41<4.06=g(\log_32),$$ the Algorithm $\bf c$
pack perfectly the squares $S^t_n$ $(n\ge4)$ into $\B$, which
completes the proof.

\medskip

\medskip

\centerline {\bf 5. DISCUSSION}

\medskip

If we increase the number of the packed squares before we start the
Algorithm $\bf c$ and do detailed analysis of the height of the
boxes, then we can decrease the constant $\log_32$. It remains an
interesting question to increase the bound $2/3$.

\begin{flushright}
{\small University of Duna\'ujv\'aros,\\
              Duna\'ujv\'aros, Hungary 2400\\
              {joosa@uniduna.hu}}
\end{flushright}
\end{document}